\documentclass[11pt]{amsart}

\usepackage{epigamath}

\usepackage[english]{babel}

\numberwithin{equation}{section}

\usepackage{enumitem}
\usepackage{tikz-cd}

\newtheorem{theorem}{Theorem}[section]

\newtheorem{proposition}[theorem]{Proposition}

\theoremstyle{definition}

\newtheorem*{question}{Question}

\theoremstyle{remark}

\def\m{\mathbb}

\EpigaVolumeYear{6}{2022} \EpigaArticleNr{10} \ReceivedOn{June 22,
2021}
\InFinalFormOn{February 4, 2022}
\AcceptedOn{March 4, 2022}

\title{Examples of surfaces with canonical map of degree 4}
\titlemark{Examples of surfaces with canonical map of degree 4}

\author{Carlos Rito}
\address{Universidade de Tr\'as-os-Montes e Alto Douro,  Quinta de Prados, 5000-801 Vila Real, Portugal}
\email{crito@utad.pt}

\authormark{C.~Rito}

\AbstractInEnglish{We give two examples of surfaces with canonical map of degree 4 onto a canonical surface.}



\acknowledgement{This research was partially supported by FCT --- Funda\c c\~ao para a Ci\^encia e a Tecnologia, I.P. ---
under the fellowship SFRH/BPD/111131/2015 and CMUP (project with reference UIDB/00144/2020).}

\begin{document}



\maketitle

\begin{prelims}

\DisplayAbstractInEnglish

\bigskip

\DisplayKeyWords

\medskip

\DisplayMSCclass







\end{prelims}


\newpage

\setcounter{tocdepth}{1}

\tableofcontents


\section{Introduction}

Let $S$ be a smooth minimal surface of general type with geometric genus $p_g\geq 3$.
Denote by ${\phi:S\dashrightarrow\m P^{p_g-1}}$ the canonical map and let $d:=\deg(\phi).$
The following result of Beauville is well-known.
\begin{theorem}[\cite{Be1}]\label{Beauv}
If the canonical image $\Sigma:=\phi(S)$ is a surface, then either:
\begin{enumerate}
\item[$(A)$] $p_g(\Sigma)=0,$ or
\item[$(B)$] $\Sigma$ is a canonical surface (i.e. it is the canonical image of a surface with birational canonical map),
in particular $p_g(\Sigma)=p_g(S)$.
\end{enumerate}
Moreover, in case $(A)$ $d\leq 36$ and in case $(B)$ $d\leq 9.$
\end{theorem}

The question of which pairs $(d,p_g)$ can actually occur has been object of study for some authors.
Several examples were given for case $(A),$ but case $(B)$ is still mysterious.
It is known that if $d>3,$ then $p_g\leq 12,$ but so far only the case $(d,p_g)=(5,4)$ has been shown to exist
(independently by Tan \cite{ta} and by Pardini \cite{Pa2}).
We refer the recent preprint by Mendes Lopes and Pardini \cite{MP} for a more detailed account on the subject.
They leave some open problems, this note is motivated by their last question.

\begin{question}
For what pairs $(d,p_g)$, with $d > 3,$ are there examples of surfaces in case $(B)$ of Theorem \ref{Beauv}?
\end{question}

Here we give examples for the cases $(d,p_g)=(4,5)$ and $(4,7)$, with canonical images a 40-nodal complete intersection
surface in $\m P^4$ and a 48-nodal complete intersection surface in $\m P^6,$ respectively
(Beauville also paid some attention to such nodal surfaces, see \cite{Be2}).

The strategy for the construction is the following. If $X$ is a surface with nodes admitting a Galois covering
$Y\rightarrow X$ ramified over the nodes and with Galois group $G$, a group with a ``big'' number of subgroups,
then we have a ``big'' number of intermediate coverings of $X.$
By computing the geometric genus $p_g$ of all involved surfaces, we may hope to find some
$\rho:W\rightarrow Z$ with $p_g(W)=p_g(Z),$ hence such that the canonical map of $W$ factors through $\rho.$

We work explicitely with the equations of a 40-nodal surface from \cite{RRS},
all computations are implemented with Magma \cite{BCP}.

\subsection*{Notation}

As usual the holomorphic Euler characteristic of a surface $S$ is denoted by $\chi(S),$ the geometric genus by $p_g(S),$
the irregularity by $q(S),$ and a canonical divisor by $K_S.$
A $(-m)$-curve is a curve isomorphic to $\m P^1$ with self-intersection $-m$.
A node of $S$ is an ordinary double point of $S.$
We say that a set of nodes of $S$ is 2-divisible if the sum $\sum A_i$ of the corresponding $(-2)$-curves
in the smooth minimal model of $S$ is 2-divisible in the Picard group.

\section{$(\mathbb Z/2)^r$-coverings}

The following result is taken from \cite[Proposition 7.6]{Ca}. See also \cite{Pa1}. 
\begin{proposition}\label{CatProp}
A normal finite $G\cong(\mathbb Z/2)^r$-covering $\pi:Y\rightarrow X$ of a smooth variety $X$ is completely
determined by the datum of
\begin{enumerate}
\item reduced effective divisors $D_{\sigma},$ for all $\sigma\in G,$ with no common components;
\item divisor linear equivalence classes $L_{\chi_1},\ldots,L_{\chi_r},$ for $\chi_1,\ldots,\chi_r$ a basis of the group of characters
$G^{\vee},$ such that $$2L_{\chi_i}\equiv\sum_{\chi_i(\sigma)=1}D_{\sigma}$$ 
\end{enumerate}
$($with additive notation for the characters$\,)$.

Conversely, given (1) and (2), one obtains a normal scheme $Y$ with a finite
$G\cong(\mathbb Z/2)^r$-covering $Y\rightarrow X,$ with branch curves the divisors $D_{\sigma}.$
\end{proposition}
The scheme $Y$ is irreducible if $\{\sigma|D_{\sigma}>0\}$ generates $G$.
We have a splitting
$$\pi_*\mathcal O_Y=\bigoplus_{\chi\in G^{\vee}}L_{\chi}^{-1}.$$

From now on, we assume that $X$ and $Y$ are surfaces. If each $D_{\sigma}$ is smooth and $\sum D_{\sigma}$ has simple normal crossings,
then $Y$ is smooth and its invariants are
\begin{equation} \label{chiformula}
\begin{split}
\chi(\mathcal O_Y)&=2^r\chi(\mathcal O_X)+\frac{1}{2}\sum_{\chi\in G^{{\vee}*}} \left(L_{\chi}^2+K_X\cdot L_{\chi}\right),\\
p_g(Y)&=p_g(X)+\sum_{\chi\in G^{{\vee}*}} h^0(X,\mathcal O_X(K_X+L_{\chi})).
\end{split}
\end{equation}
Let $R_{\sigma}$ be the support of $\pi^*(D_{\sigma})$.
The Hurwitz formula gives $$K_Y\equiv\pi^*(K_X)+\sum_{\sigma\in G^*}R_{\sigma}.$$

Now assume that the $D_{\sigma}$ are disjoint $(-2)$-curves. Then the $R_{\sigma}$ are disjoint $(-1)$-curves,
the canonical map of $Y$ factors through the covering $Y\rightarrow X$ if and only if $p_g(Y)=p_g(X),$
and one has a commutative diagram
\begin{center}
\begin{tikzcd}[arrows={-Stealth}]
  Y\rar\dar & Y'\dar \\%
  X\rar & X'
\end{tikzcd}
\end{center}
where $Y\rightarrow Y'$ is the contraction of the $(-1)$-curves $R_{\sigma},$
the surface $X'$ has nodes corresponding to the $(-2)$-curves of $X,$
and $Y'\rightarrow X'$ is a $(\mathbb Z/2)^r$-covering ramified on those nodes.
In this case Equation (\ref{chiformula}) becomes
\begin{equation}\label{formula}
\chi(\mathcal O_Y)=2^r\left(\chi(\mathcal O_X)-m/8\right)
\end{equation}
where $m$ is the number of nodes of $X'.$

\section{Construction}

Let $X_{40}$ be the surface in $\mathbb P^4$ given by the equations
\begin{equation} \label{eqs}
\begin{aligned}
5\left(x^2+y^2+z^2+w^2+t^2\right)-7\left(x+y+z+w+t\right)^2&=0 \\
4\left(x^4+y^4+z^4+w^4+t^4+h^4\right)-\left(x^2+y^2+z^2+w^2+t^2+h^2\right)^2&=0
\end{aligned}
\end{equation}
where
$$h:=-(x+y+z+w+t).$$
It is the canonical model of a surface with invariants $p_g=5,$ $q=0$ and $K^2=8.$
The above quartic $I$ is classically known as the Igusa quartic; its singular set is the union of 15 lines.
The quadric meets these lines transversally, and is tangent to $I$ at 10 smooth points,
thus the singular set of $X_{40}$ is the union of 40 nodes $N_1,\ldots,N_{40}$ (for more details see \cite{RRS}).

Let $\widetilde X_{40}$ be the smooth minimal model of $X_{40}$ and denote by
$A_i$ the $(-2)$-curves in $\widetilde X_{40}$ corresponding to the nodes $N_i,$ $i=1,\ldots,40.$
Let $a,b,c$ be the canonical generators of the group $(\m Z/2)^3$ and, for $i,j,k\in \mathbb Z/2,$
let $\chi_{ijk}$ denote the character which takes the value $i,j,k$ on $a,b,c,$ respectively.
We show in Section \ref{computations} that one can write
$$A_1+\cdots+A_{40}=D_a+D_b+D_c+D_{abc}+D_{bc}+D_{ac}+D_{ab}$$
where each of $D_a,D_b,D_c,D_{abc}$ is a sum of $4$ $(-2)$-curves,
each of $D_{bc},D_{ac},D_{ab}$ is a sum of $8$ $(-2)$-curves,
and such that there exist divisors $L_{100},L_{010},L_{001}$ satisfying:
\begin{equation}\label{divs1}
\begin{aligned}
D_a+D_{abc}+D_{ac}+D_{ab}&\equiv 2L_{100} \\
D_b+D_{abc}+D_{bc}+D_{ab}&\equiv 2L_{010} \\
D_c+D_{abc}+D_{bc}+D_{ac}&\equiv 2L_{001}.
\end{aligned}
\end{equation}
It follows from Proposition~\ref{CatProp} that these data
define a $(\m Z/2)^3$-covering $\pi:\widetilde Y\rightarrow\widetilde X_{40}$ branched on the $(-2)$-curves
$A_i,$ equivalently a $(\m Z/2)^3$-covering $\psi:Y\rightarrow X_{40}$ branched on the nodes of $X_{40}$
(the surface $Y$ is minimal because $X_{40}$ is minimal and $\psi$ is \'etale in codimension 1).
In particular there exist divisors $L_{111},L_{110},L_{101},L_{011}$ such that:
\begin{equation}\label{divs2}
\begin{aligned}
D_a+D_b+D_c+D_{abc}&\equiv 2L_{111}\\
D_a+D_b+D_{bc}+D_{ac}&\equiv 2L_{110}\\
D_a+D_c+D_{bc}+D_{ab}&\equiv 2L_{101}\\
D_b+D_c+D_{ac}+D_{ab}&\equiv 2L_{011}.
\end{aligned}
\end{equation}
One has
\begin{equation*}
2L_{ijk}\equiv\sum_{\chi_{ijk}(\sigma)=1}D_{\sigma}.
\end{equation*}

Since $\psi$ is ramified only on nodes, we have $K_Y\equiv\psi^*(K_{X_{40}})$ and then $K_Y^2=8K_{X_{40}}^2=64.$
We show in Section \ref{computations} that
$$h^0\left({\widetilde X_{40}},\mathcal O_{\widetilde X_{40}}\left(K_{\widetilde X_{40}}+L_{111}\right)\right)=2$$
and
$$h^0\left({\widetilde X_{40}},\mathcal O_{\widetilde X_{40}}\left(K_{\widetilde X_{40}}+L_{ijk}\right)\right)=0\ \ {\rm for}\ \ {ijk}\ne 111,$$
thus $$p_g(Y)=p_g(X_{40})+2+0+\cdots+0=7.$$

We get from (\ref{formula}) that $\chi(Y)=8(6-5)=8,$ thus $q(Y)=0.$

The covering $\psi$ factors as

\begin{center}
\begin{tikzcd}[arrows={-Stealth}]
  Y\rar\dar & Y_{32}\rar & Y_{48}\dar \\%
  X_{16}\rar & X_{32}\rar & X_{40}
\end{tikzcd}
\end{center}
with $Y_{48}$ and $X_{16}$ given by the quotients by the groups $\langle ab,ac\rangle$ and $\langle c\rangle,$ respectively
(the subscript $n$ means a surface with singular set the union of $n$ nodes).
All these surfaces are regular because $q(Y)=0.$

It follows from (\ref{formula}) that $\chi(X_{16})=4(6-36/8)=6,$ thus $p_g(X_{16})=p_g(X_{40})=5,$
and we conclude that
\begin{center}
the $(\m Z/2)^2$-covering $X_{16}\rightarrow X_{40}$ is the canonical map of $X_{16}.$
\end{center}
Analogously, $p_g(Y)=p_g(Y_{48})=7$ and we claim that
\begin{center}
the $(\m Z/2)^2$-covering $Y\rightarrow Y_{48}$ is the canonical map of $Y.$
\end{center}
For this it suffices to show that $Y_{48}$ is a canonical surface.

Since the canonical system of $Y_{48}$ contains the pullback of the canonical system of $X_{40}$
and since $p_g(Y_{48})>p_g(X_{40}),$ the canonical map of $Y_{48}$ must be birational.
But we can be more precise.
We follow Beauville \cite{Be2} and show that $Y_{48}$ can be embedded in $\m P^6$ as a complete
intersection of 4 quadrics in the following way. The linear system $L$ of quadrics through the branch locus of the covering
$Y_{48}\rightarrow X_{40}$ (16 nodes) is of dimension 2.
Using computer algebra it is not difficult to show that $L$ contains quadrics $B,C,D$ such that
the surface $X_{40}$ is given by $Q=0,$ $B^2-CD=0,$ where $Q$ is the quadric from (\ref{eqs})
(we write the quadrics as general elements of $L$, thus depending on some parameters;
then we obtain a variety on these parameters by imposing that the hypersurfaces $Q=0$ and $B^2-CD=0$
are tangent at the 24 nodes of $X_{40}$ which are disjoint from the 16 nodes of $B^2-CD=0$;
finally we compute points in this variety).

Then $Y_{48}$ is given in $\m P^6(x,y,z,w,t,u,v)$ by equations
$$u^2-C=v^2-D=uv-B=Q=0.$$
We give these equations in Section \ref{Y48} and verify that $Y_{48}$ is as stated.

Let us explain how we find 2-divisible sets of nodes in $X_{40}.$
The surface $X_{40}$ contains 40 {\em tropes}, which are hyperplane sections $H_i=2T_i$
with $T_i\subset X_{40}$ a reduced curve through 12 nodes of $X_{40},$ and smooth at these points.
Thus in $\widetilde X_{40}$ the pullback of such a trope can be written as
$$\widetilde H_i=2\widehat T_i+\sum_{j\in J}A_j\, ,\ \text{with}\ \# J=12.$$
Thus for each pair of tropes the sum of nodes contained in their union
and not contained in their intersection is 2-divisible.

Using these 2-divisibilities, the strategy for finding configurations as in (\ref{divs1}) is simple:
we have used a computer algorithm to list and check possibilities.

\section{Computations}

The computations below are implemented with Magma V2.26-5.

\subsection{The covering $Y\rightarrow X_{40}$}\label{computations}

We start by defining the surface $X_{40}$ and its singular set.

\begin{verbatim}
K:=Rationals();
R<r>:=PolynomialRing(K);
K<r>:=ext<K|r^2 + 15>;
P<x,y,z,w,t>:=ProjectiveSpace(K,4);
h:=-x-y-z-w-t;
Q:=5*(x^2+y^2+z^2+w^2+t^2)-7*(x+y+z+w+t)^2;
I:=4*(x^4+y^4+z^4+w^4+t^4+h^4)-(x^2+y^2+z^2+w^2+t^2+h^2)^2;
X40:=Surface(P,[Q,I]);
SX40:=SingularSubscheme(X40);
\end{verbatim}

The partition of the 40 nodes:

\begin{verbatim}
Da:={P![3,3,-2,-2,3],P![4,-r+1,r-5,-r+1,4],
     P![-r+1,4,r-5,-r+1,4],P![r+1,r+1,-r-5,4,4]};
Db:={P![2,-3,-3,-3,2],P![4,r+1,r+1,-r-5,4],
     P![-r-5,r-5,r-5,-r-5,10],P![r-5,-r+1,-r+1,4,4]};
Dc:={P![-3,-3,2,-3,2],P![-r+1,-r+1,r-5,4,4],
     P![r-5,r-5,-r-5,-r-5,10],P![r+1,r+1,4,-r-5,4]};
Dabc:={P![-2,3,3,-2,3],P![-r-5,r+1,r+1,4,4],
       P![r-5,4,-r+1,-r+1,4],P![r-5,-r+1,4,-r+1,4]};
Dbc:={P![-2,-2,3,3,3],P![3,-2,-2,3,3],
      P![4,-r-5,r+1,r+1,4],P![4,-r+1,-r+1,r-5,4],
      P![4,r+1,-r-5,r+1,4],P![-r-5,r+1,4,r+1,4],
      P![-r+1,-r+1,4,r-5,4],P![r+1,-r-5,4,r+1,4]};
Dac:={P![-3,2,-3,-3,2],P![3,-2,3,-2,3],
      P![4,r-5,-r+1,-r+1,4],P![-r+1,r-5,4,-r+1,4],
      P![-r+1,r-5,-r+1,4,4],P![r-5,-r-5,r-5,-r-5,10],
      P![r+1,4,r+1,-r-5,4],P![r+1,-r-5,r+1,4,4]};
Dab:={P![-3,-3,-3,2,2],P![-2,3,-2,3,3],
      P![-r-5,4,r+1,r+1,4],P![-r+1,4,-r+1,r-5,4],
      P![-r-5,-r-5,r-5,r-5,10],P![-r-5,r-5,-r-5,r-5,10],
      P![r-5,-r-5,-r-5,r-5,10],P![r+1,4,-r-5,r+1,4]};
\end{verbatim}

Verification that these are in fact the nodes:

\begin{verbatim}
&join[Da,Db,Dc,Dabc,Dbc,Dac,Dab] eq SingularPoints(X40);
HasSingularPointsOverExtension(X40) eq false;
\end{verbatim}

Some of the tropes of $X_{40}:$

\begin{verbatim}
tropes:=[
    6*x + (-r - 9)*y + (r - 9)*z + (r - 9)*w + (-r - 9)*t,
    16*x + (-r - 9)*y + 16*z + (3*r + 11)*w + (3*r + 11)*t,
    16*x + (r - 9)*y + 16*z + (-3*r + 11)*w + (-3*r + 11)*t,
    6*x + (r - 9)*y + (-r - 9)*z + (r - 9)*w + (-r - 9)*t,
    16*x + (3*r + 11)*y + 16*z + (3*r + 11)*w + (-r - 9)*t,
    16*x + (-3*r + 11)*y + (-3*r + 11)*z + (r - 9)*w + 16*t,
    x + y + w,
    16*x + (r - 9)*y + (-3*r + 11)*z + (-3*r + 11)*w + 16*t,
    x + z + w
];
\end{verbatim}

The reduced subscheme of these tropes:

\begin{verbatim}
red:=[ReducedSubscheme(Scheme(X40,q)):q in tropes];
&and[Degree(q) eq 4:q in red];
\end{verbatim}

They are smooth at the nodes of $X_{40}$:

\begin{verbatim}
&and[Dimension(SingularSubscheme(q) meet SX40) eq -1:q in red];
\end{verbatim}

Two 2-divisible disjoint sets of 20 nodes, which confirm that the 40 nodes are 2-divisible:

\begin{verbatim}
s1:=Points(Scheme(SX40,tropes[1]*tropes[2])) diff
    Points(Scheme(SX40,[tropes[1],tropes[2]]));
s2:=Points(Scheme(SX40,tropes[6]*tropes[7])) diff
    Points(Scheme(SX40,[tropes[6],tropes[7]]));
&and[#s1 eq 20,#s2 eq 20,#(s1 join s2) eq 40];
\end{verbatim}

We compute three 2-divisible sets of 24 nodes:

\begin{verbatim}
Sets:=[];
for q in [[2,5],[1,4],[3,8]] do
  pts:=Points(Scheme(SX40,tropes[q[1]]*tropes[q[2]])) diff
       Points(Scheme(SX40,[tropes[q[1]],tropes[q[2]]]));
  Append(~Sets,SingularPoints(X40) diff pts);
end for;
\end{verbatim}

and use these sets to check the divisibilities in (\ref{divs1}):

\begin{verbatim}
Da join Dabc join Dac join Dab eq Sets[1];
Db join Dabc join Dbc join Dab eq Sets[2];
Dc join Dabc join Dbc join Dac eq Sets[3];
\end{verbatim}

Now we show that 
\[h^0\left(\widetilde X_{40},\mathcal O_{\widetilde X_{40}}\left(K_{\widetilde X_{40}}+L_{111}\right)\right)=2.\]
Let $N_1,\ldots,N_{16}$ be the nodes in $D_a+D_b+D_c+D_{abc}$ and $A_1,\ldots,A_{16}$ be the corresponding $(-2)$-curves.
Let $H_1,H_2$ be the tropes whose pullback to $\widetilde X_{40}$ is
$$\widetilde H_1+\widetilde H_2=2\widehat T_1+2\widehat T_2+\sum_{i=1}^{16}A_i+2\sum_{i=17}^{20}A_i,$$
with $A_{17},\ldots,A_{20}\in \widetilde H_1\cap\widetilde H_2.$
Then
$$\sum_{i=1}^{16}A_i\equiv 2L_{111},\ \ \ {\rm with}\ \ \ K_{\widetilde X_{40}}+L_{111}\equiv 2\widetilde H-\widehat T_1-\widehat T_2-\sum_{i=17}^{20}A_i.$$
We compute below that the system of quadrics through the curves $T_1,T_2\subset\m P^4$ is generated by 2 elements, modulo the quadric $Q$.
For $i=17,\ldots 20,$ the fact $\left(2\widetilde H-\widehat T_1-\widehat T_2\right)\cdot A_i<0$ implies that $A_i$ is contained in the base component of the
linear system $\left|2\widetilde H-\widehat T_1-\widehat T_2\right|.$
This gives
$h^0\left(\widetilde X_{40},\mathcal O_{\widetilde X_{40}}\left(K_{\widetilde X_{40}}+L_{111}\right)\right)=2.$

\begin{verbatim}
T1:=ReducedSubscheme(Scheme(X40,tropes[2]));
T2:=ReducedSubscheme(Scheme(X40,tropes[9]));
pts:=Points(SX40 meet (T1 join T2)) diff
     Points(SX40 meet T1 meet T2);
pts eq (Da join Db join Dc join Dabc);
L:=LinearSystem(LinearSystem(P,2),T1 join T2);
#Sections(LinearSystemTrace(L,X40)) eq 2;
\end{verbatim}

Let us show that
\[h^0\left(\widetilde X_{40},\mathcal O_{\widetilde X_{40}}\left(K_{\widetilde X_{40}}+L_{ijk}\right)\right)=0\]
for $ijk\ne 111$. Suppose the opposite. Let $A_1,\ldots,A_{24}$ be the corresponding $(-2)$-curves.
Then there is a curve $E\in\left|K_{\widetilde X_{40}}+L_{ijk}\right|,$
and $E\cdot A_i=-1$ implies that the linear system
$\left|K_{\widetilde X_{40}}+L_{ijk}-\sum_{j=1}^{24}A_j\right|=\left|K_{\widetilde X_{40}}-L_{ijk}\right|$
is nonempty.
Therefore $\left|2K_{\widetilde X_{40}}-\sum_{j=1}^{24}A_j\right|$ is nonempty, which implies that
there is at least one quadric in $\m P^4$ through the corresponding nodes $N_1,\ldots,N_{24}$
(modulo the quadric $Q$).
We show below that this does not happen.

\begin{verbatim}
Sets:=[
Da join Dabc join Dac join Dab,
Db join Dabc join Dbc join Dab,
Dc join Dabc join Dbc join Dac,
Da join Db join Dbc join Dac,
Da join Dc join Dbc join Dab,
Db join Dc join Dac join Dab
];
for q in Sets do
  L:=LinearSystem(LinearSystem(P,2),[P!x:x in q]);
  #Sections(LinearSystemTrace(L,X40)) eq 0;
end for;
\end{verbatim}

\subsection{The surface $Y_{48}$}\label{Y48}

Here we give the equations of $Y_{48}$ as a complete intersection of 4 quadrics in $\m P^6.$
We start by defining $\m P^6$ over a certain number field.

\begin{verbatim}
K:=Rationals(); R<x>:=PolynomialRing(K);
K<r,m>:=ext<K|x^2 + 15,x^2 - 95/42*x + 2855/2646>;
R<n>:=PolynomialRing(K);
K<n>:=ext<K|
n^2 + 443889677/206391214080000*r - 46942774543/619173642240000>;
P6<x,y,z,w,t,u,v>:=ProjectiveSpace(K,6);
\end{verbatim}

The three quadrics $B,C,D$:

\begin{small}
\begin{verbatim}
B:=(675/4802*r+334125/33614)*n*x*z+(-389475/67228*r+3266325/67228)*n*x*w+
 (34425/9604*r+451575/67228)*n*y*w+(-389475/67228*r+3266325/67228)*n*z*w+
 (-62100/16807*r+348300/16807)*n*w^2+(239625/33614*r+1541025/33614)*n*x*t
 +(-8100/2401*r+137700/16807)*n*y*t+(239625/33614*r+1541025/33614)*n*z*t
 +(6075/9604*r+3007125/67228)*n*w*t+(71550/16807*r+319950/16807)*n*t^2;
C:=x*y+1/154*(126*m-181)*y^2+1/42*(-42*m+95)*x*z+y*z+(1/1540*(14*m-25)*r
 +1/924*(-798*m+1997))*x*w+(1/420*(42*m-65)*r+1/308*(-294*m+767))*y*w
 +(1/1540*(14*m-25)*r+1/924*(-798*m+1997))*z*w+(1/385*(-119*m+185)*r
 +1/462*(-168*m+311))*w^2+(1/1540*(-14*m+25)*r+1/924*(-798*m+
 1997))*x*t+(1/420*(-42*m+65)*r+1/308*(-294*m+767))*y*t+(1/1540*(-14*m
 +25)*r+1/924*(-798*m+1997))*z*t+1/154*(126*m-71)*w*t+(1/385*(119*m-
 185)*r+1/462*(-168*m+311))*t^2;
D:=x*y+1/77*(-63*m+52)*y^2+m*x*z+y*z+(1/2310*(-21*m+10)*r+1/154*(133*m+
 32))*x*w+(1/70*(-7*m+5)*r+1/154*(147*m+51))*y*w+(1/2310*(-21*m+
 10)*r+1/154*(133*m+32))*z*w+(1/2310*(714*m-505)*r+1/154*(56*m-
 23))*w^2+(1/2310*(21*m-10)*r+1/154*(133*m+32))*x*t+(1/70*(7*m-5)*r
 +1/154*(147*m+51))*y*t+(1/2310*(21*m-10)*r+1/154*(133*m+32))*z*t+
 1/77*(-63*m+107)*w*t+(1/2310*(-714*m+505)*r+1/154*(56*m-23))*t^2;
\end{verbatim}
\end{small}

We obtain alternative equations for $X_{40}$:

\begin{verbatim}
F:=B^2-C*D;
Q:=5*(x^2+y^2+z^2+w^2+t^2)-7*(x+y+z+w+t)^2;
X:=Scheme(P6,[F,Q,u,v]);
h:=-x-y-z-w-t;
I:=4*(x^4+y^4+z^4+w^4+t^4+h^4)-(x^2+y^2+z^2+w^2+t^2+h^2)^2;
X40:=Scheme(P6,[Q,I,u,v]);
X eq X40;
\end{verbatim}

And finally the equations of $Y_{48}$ in $\m P^6$:

\begin{verbatim}
Y48:=Surface(P6,[u^2-C,v^2-D,u*v-B,Q]);
SY48:=SingularSubscheme(Y48);
Dimension(SY48) eq 0;
Degree(SY48) eq 48;
Degree(ReducedSubscheme(SY48)) eq 48;
\end{verbatim}


\end{document}